\documentclass[10pt]{amsart}
\input{epsf}
\usepackage{amssymb,latexsym}
\topskip  \numberwithin{equation}{section}
\newtheorem{theorem}{Theorem}[section]

\newtheorem{remark}[theorem]{Remark}

\newtheorem{example}[theorem]{Example}

\begin{document}

\pagenumbering{arabic}
\pagestyle{headings}
\def\sof{\hfill\rule{2mm}{2mm}}
\def\llim{\lim_{n\rightarrow\infty}}

\title {Counting the occurrences of generalized patterns in words generated by a morphism}
\maketitle

\begin{center}Sergey Kitaev and Toufik Mansour \footnote{Research
financed by EC's
IHRP Programme, within the Research Training Network "Algebraic
Combinatorics in Europe", grant HPRN-CT-2001-00272}
\end{center}
\begin{center}
Matematik, Chalmers tekniska h\"ogskola och G\"oteborgs universitet,\\
S-412~96 G\"oteborg, Sweden

{\tt kitaev@math.chalmers.se, toufik@math.chalmers.se}
\end{center}

\def\mn{\mbox{-}}
\def\SS{{\mathcal S}}

\section*{Abstract}
We count the number of occurrences of certain patterns
in given words. We choose these words to be the set of all finite
approximations of a sequence generated by a morphism with certain
restrictions. The patterns in our considerations are either
classical patterns $1\mn 2$, $2 \mn 1$, $1\mn 1\mn \cdots \mn 1$,
or arbitrary generalized patterns without internal dashes, in
which repetitions of letters are allowed. In particular, we find
the number of occurrences of the patterns $1\mn 2$, $2\mn 1$,
$12$, $21$, $123$ and $1\mn 1\mn \cdots \mn 1$ in the words
obtained by iterations of the morphism $1\rightarrow 123$,
$2\rightarrow 13$, $3\rightarrow 2$, which is a classical example
of a morphism generating a nonrepetitive sequence.
\thispagestyle{empty} 


\section{Introduction and Background}\label{introduction}

We write permutations as words $\pi=a_1 a_2\cdots a_n$, whose letters are distinct and usually consist of the integers $1,2,\ldots,n$.

An occurrence of a pattern $p$ in a permutation $\pi$ is ``classically'' defined as a subsequence in $\pi$ (of the same length as the length of $p$) whose letters are in the same relative order as those in $p$. Formally speaking, for $r \leq n$, we say that a permutation $\sigma$ in the symmetric group ${\mathcal S}_n$ has an occurrence of the pattern $p \in {\mathcal S}_r$ if there exist $1 \leq i_1 < i_2 < \cdots < i_r \leq n$ such that $p = \sigma(i_1)\sigma(i_2) \ldots \sigma(i_r)$ in reduced form. The {\em reduced form} of a permutation $\sigma$ on a set $\{ j_1,j_2, \ldots ,j_r \}$, where $j_1 < j_2 < \cdots <j_r$, is a permutation ${\sigma}_1$ obtained by renaming the letters of the permutation $\sigma$ so that $j_i$ is renamed $i$ for all $i \in \{1, \ldots ,r\}$. For example, the reduced form of the permutation 3651 is 2431. The first case of classical patterns studied was that of permutations avoiding a pattern of length 3 in $\SS_3$. Knuth \cite{Knuth} found that, for any $\tau\in\SS_3$, the number $|\SS_n(\tau)|$ of $n$-permutations avoiding $\tau$ is $C_n$, the $n$th Catalan number.
Later, Simion and Schmidt \cite{SimSch} determined the number $|\SS_n(P)|$ of permutations in $\SS_n$ simultaneously avoiding any given set of patterns $P\subseteq\SS_3$.

In \cite{BabStein} Babson and Steingr\'{\i}msson introduced {\em generalised permutation patterns} that allow the requirement that two adjacent letters in a pattern must be adjacent in the permutation. In order to avoid confusion we write a "classical" pattern, say $231$, as $2$-$3$-$1$, and if we write, say $2$-$31$, then we mean that if this pattern occurs in the permutation, then the letters in the permutation that correspond to $3$ and $1$ are adjacent. For example, the permutation $\pi=516423$ has only one occurrence of the pattern $2$-$31$, namely the subword 564, whereas the pattern $2$-$3$-$1$ occurs, in addition, in the subwords 562 and 563. A motivation for introducing these patterns in \cite{BabStein} was the study of Mahonian statistics. A number of interesting results on generalised patterns were obtained in \cite{Claes}. Relations to several well studied combinatorial structures, such as set partitions, Dyck paths, Motzkin paths and involutions, were shown there.

Burstein \cite{Burstein} considered words instead of permutations. In particular, he found the number $|[k]^n(P)|$ of words of length $n$ in $k$-letter alphabet that avoid each pattern from a set $P\subseteq\SS_3$ simultaneously. Burstein and Mansour \cite{BurMans1} (resp. \cite{BurMans2,BurMans3}) considered forbidden patterns (resp. generalized patterns) with repeated letters.

The most attention, in the papers on classical or generalized patterns, is paid to counting exact formulas and/or generating functions for the number of words or permutations avoiding, or having $k$ occurrences of, certain pattern. In this paper we suggest another problem, namely counting the number of occurrences of a particular pattern $\tau$ in given words. We choose these words to be a set of all finite approximations (to be defined below) of a sequence generated by a morphism with certain restrictions. A motivation for such a choice is big interest in studying classes of sequences and words that are defined by iterative schemes~\cite{Lothaire,Salomaa}. The pattern $\tau$ in our considerations is either a classical pattern from the set $\{ 1\mn 2, 2 \mn 1, 1\mn 1\mn \cdots \mn 1 \}$, or an arbitrary generalized pattern without internal dashes, in which repetitions of letters are allowed. In particular, we find that there are $(3\cdot 4^{n-1}+2^n)$ occurrences of the pattern $1\mn 2$ in the $n$-th finite approximation of the sequence $w$ defined below, which is a classical example of a nonrepetitive sequence.

Let $\Sigma$ be an alphabet and ${\Sigma}^{\star}$ be the set of all words of $\Sigma$. A map $\varphi: {\Sigma}^{\star} \rightarrow {\Sigma}^{\star}$ is called a {\em morphism}, if we have $\varphi(uv)=\varphi(u)\varphi(v)$ for any $u,v \in {\Sigma}^{\star}$. It is easy to see that a morphism $\varphi$ can be defined by defining $\varphi(i)$ for each $i \in \Sigma$. The set of all rules $i \rightarrow \varphi(i)$ is called a {\em substitution system}. We create words by starting with a letter from the alphabet $\Sigma$ and iterating the substitution system. Such a substitution system is called a {\em D0L (Deterministic, with no context Lindenmayer) system} \cite{LindRoz}. D0L systems are classical objects of formal language theory. They are interesting from mathematical point of view~\cite{Frid}, but also have applications in theoretical biology~\cite{Lind}. Let $|X|$ denote the length of a word $X$, that is the number of letters in~$X$.

Suppose a word $\varphi(a)$ begins with $a$ for some $a \in \Sigma$, and that the length of ${\varphi^k}(a)$ increases without bound. The symbolic sequence $\lim\limits_{ k \to \infty }{\varphi^k}(a)$ is said to be {\em generated} by the morphism $\varphi$. In particular, $\lim\limits_{ k \to \infty }{\varphi^k}(a)$ is a {\em fixed point} of $\varphi$. However, in this paper we are only interesting in the {\em finite approximations} of $\lim\limits_{ k \to \infty }{\varphi^k}(a)$, that is in the words ${\varphi^k}(a)$ for $k=1,2,\ldots$.

An example of a sequence generated by a morphism can be the following sequence~$w$. We create words by starting with the letter 1 and iterating the substitution system $\phi_{w}$: $1\rightarrow 123$, $2\rightarrow 13$, $3\rightarrow 2$. Thus, the initial letters of $w$ are 123132123213.... This sequence was constructed in connection with the problem of constructing a nonrepetitive sequence on a 3-letter alphabet, that is, a sequence that does not contain any subwords of the type $XX=X^2$, where $X$ is any non-empty word over a 3-letter alphabet. The sequence $w$ has that property. The question of the existence of such a sequence, as well as the questions of the existence of sequences avoiding other kinds of repetitions, were studied in algebra \cite{Adian,Justin,Kol}, discrete analysis \cite{Carpi,Dekk,Evdok,Ker,Pleas} and in dynamical systems \cite{MorseHedl}. In Examples~\ref{ex1},~\ref{ex4} and~\ref{ex5} we give the number of occurrences of the patterns $1\mn 2$, $2\mn 1$, $1\mn 1\mn \cdots \mn 1$, $12$, $123$ and $21$ in the finite approximations of~$w$.

To proceed further, we need the following definitions. Let
$N^{\tau}_{\phi}(n)$ denote the number of occurrences of the
pattern $\tau$ in a word generated by some morphism $\phi$ after
$n$ iterations. We say that an occurrence of $\tau$ is {\em
external} for a pair of words $(X,Y)$, if this occurrence starts
in $X$ and ends in $Y$. Also, an occurrence of $\tau$ for a word
$X$ is {\em internal}, if this occurrence starts and ends in this
$X$.


\section{Patterns 1-2, 2-1 and 1-1-...-1}

\begin{theorem}\label{generalTheorem}
Let $\mathcal{A}$ $=\{ 1 ,2,\ldots, k \}$ be an alphabet, where $k\geq 2$ and a pattern $\tau\in \{ 1\mn 2,2\mn 1 \}$. Let $X_1$ begins with the letter 1 and consists of $\ell$ copies of each letter $i\in \mathcal{A}$ {\rm(}$\ell \geq 1${\rm)}. Let a morphism $\phi$ be such that
$$1 \rightarrow X_1, \ 2 \rightarrow X_2, \ 3 \rightarrow X_3, \ldots , k \rightarrow X_k,$$
where we allow $X_i$ to be the empty word $\epsilon$ for $i=2,3,\ldots, k$ {\rm(}that is, $\phi$ may be an erasing morphism{\rm)}, $\displaystyle\sum_{i=2}^{k}|X_i|=k\cdot d$, and each letter from $\mathcal{A}$ appears in the word $X_2X_3\ldots X_k$ exactly $d$ times. Besides, let $e_{i,j}$ {\rm(}resp. $e_i${\rm)} be the number of external occurrences of $\tau$ for $(X_i,X_j)$ {\rm(}resp. $(X_i,X_i)${\rm)}, where $i\neq j$. Let $s_i$ be the number of internal occurrences of $\tau$ in $X_i$. In particular, $s_i=e_i=e_{i,j}=e_{j,i}=0$, whenever $X_i=\epsilon$; also, $e_i=|X_i|\cdot (|X_i|-1)/2$, whenever there are no repetitive letters in $X_i$. Then $N^{\tau}_{\phi}(1)=s_1$ and for $n\geq 2$, $N^{\tau}_{\phi}(n)$ is given by
$$(d+\ell)^{n-2}\sum_{i=1}^{k}s_i+{(d+\ell)^{n-2} \choose 2}\sum_{i=1}^{k}e_i+(d+\ell)^{2n-4}\sum_{1\leq i < j \leq k}e_{i,j}.$$
\end{theorem}

\begin{proof}
We assume that $\tau=1\mn 2$. All the considerations for this $\tau$ remain the same for the case $\tau=2\mn 1$.

If $n=1$ then the statement is trivial.

Suppose $n\geq2$. Using the fact that $X_1X_2X_3\ldots X_k$, has exactly $d+\ell$ occurrences of each letter $i$, $i=1,2,\ldots,k$, one can prove by induction on $n$, that the word $\phi^n(1)$ is a permutation of $(d+\ell)^{n-2}$ copies of each word $X_i$, where $i=1,2,\ldots,k$. This implies, in particular, that $|\phi^n(1)|=k\cdot(d+\ell)^{n-1}$.

An occurrence of $\tau$ in $\phi^n(1)$ can be either internal, that is when $\tau$ occurs inside a word $X_i$, or external, which means that $\tau$ begins in a word $X_i$ and ends in another word $X_j$. In the first of these cases, since there are $(d+\ell)^{n-2}$ copies of each $X_i$, we have $(d+\ell)^{n-2}\sum_{i=1}^{k}s_i$ possibilities. In the second case, either $i=j$, which gives ${(d+\ell)^{n-2}\choose 2}\sum_{i=1}^{k}e_i$ possibilities, or $i\neq j$, in which case there are $(d+\ell)^{n-2}$ possibilities to choose $X_i$ (resp. $X_j$) among $(d+\ell)^{n-2}$ copies of $X_i$ (resp. $X_j$), and using the fact that $e_{i,j}=e_{j,i}$ (the order in which the words $X_i$ and $X_j$ occur in $\phi^n(1)$ is unimportant), we have $(d+\ell)^{2n-4}\sum_{1\leq i < j \leq k}e_{i,j}$ possibilities. Summing all the possibilities, we finish the proof.
\end{proof}

Let $s$ (resp. $e$) denote the vector $(s_1, s_2, \ldots, s_k)$ (resp. $(e_1, e_2, \ldots, e_k)$), where $s_i$ and $e_j$ are defined in Theorem~\ref{generalTheorem}. All of the following examples are corollaries to Theorem~\ref{generalTheorem}.

\begin{example}\label{ex1} If we consider the morphism $\phi_w$ defined in Section~\ref{introduction} and the pattern $\tau=1\mn 2$ then $d=\ell=1$, $s=(3,1,0)$, $e=(3,1,0)$ and $e_{1,2}=e_{2,1}=2$, $e_{1,3}=e_{3,1}=1$, $e_{2,3}=e_{3,2}=1$. Hence, the number of occurrences of $\tau$ is given by $N^{1\mn 2}_{\phi_w}(1)=3$ and, for $n\geq 2$, $N^{1\mn 2}_{\phi_w}(n)=(3\cdot 4^{n-1} + 2^n)/2$. If $\tau=2\mn 1$ then $s=(0,0,0)$, $e=(3,1,0)$ and $e_{1,2}=e_{2,1}=2$, $e_{1,3}=e_{3,1}=1$, $e_{2,3}=e_{3,2}=1$. Hence, $N^{2\mn 1}_{\phi_w}(1)=0$ and, for $n\geq 2$, $N^{2\mn 1}_{\phi_w}(n)=(3\cdot 4^{n-1} - 2^n)/2$.
\end{example}

\begin{example} If we consider the morphism $\phi$: $1\rightarrow 1324$, $2\rightarrow \epsilon$, $3\rightarrow 14$, and $4\rightarrow 23$ then for the pattern $\tau=1\mn 2$, we have $d=\ell=1$, $s=(5,0,1,1)$, $e=(6,0,1,1)$, and $e_{i,j}$, for $i\neq j$, are elements of the matrix
$$
\left(
\begin{array}{cccc}
- & 0 & 3 & 3 \\
0 & - & 0 & 0 \\
3 & 0 & - & 2 \\
3 & 0 & 2 & -
\end{array}
\right).
$$
Hence, $N^{1\mn 2}_{\phi}(1)=5$ and, for $n\geq 2$, $N^{1\mn 2}_{\phi}(n)=3\cdot 4^{n-1}+11\cdot 2^{n-2}$.
\end{example}

\begin{example} If we consider the morphism $\phi$: $1\rightarrow 13542$, $2\rightarrow 423$, $3\rightarrow \epsilon$, $4\rightarrow 5115$, and $5\rightarrow 234$ then for the pattern $\tau=1\mn 2$, we have $\ell=1$, $d=2$, $s=(6,1,0,2,3)$, $e=(10,3,0,4,3)$, and $e_{i,j}$, for $i\neq j$, are elements of the matrix
$$
\left(
\begin{array}{ccccc}
- & 6 & 0 & 8 & 6 \\
6 & - & 0 & 6 & 3 \\
0 & 0 & - & 0 & 0 \\
8 & 6 & 0 & - & 6 \\
6 & 3 & 0 & 6 & -
\end{array}
\right).
$$
Hence, $N^{1\mn 2}_{\phi}(1)=6$ and, for $n\geq 2$, $N^{1\mn 2}_{\phi}(n)=5\cdot 9^{n-1} + 2\cdot 3^{n-2}$.
\end{example}

Using the proof of Theorem~\ref{generalTheorem}, we have the following.

\begin{theorem}\label{multipattern}
Let a morphism $\phi$ satisfy all the conditions in the statement of Theorem~\ref{generalTheorem} and the pattern $\tau=\underbrace{1\mn 1\mn \cdots \mn 1}_{r\mbox{ times}}$. Then, for $n\geq 2$, the number of occurrences of $\tau$ in $\phi^n(1)$ is given by $k\cdot {(d+\ell)^{n-1} \choose r}$, whereas for $n=1$, by $k\cdot{\ell \choose r}$.
\end{theorem}

\begin{proof}
From the proof of Theorem~\ref{generalTheorem}, we have that if $n\geq 2$ (resp. $n=1$) then $\phi^n(1)$ has exactly $(d+\ell)^{n-1}$ (resp. $\ell$) copies of each letter from ${\mathcal A}$. We can choose $r$ of them in ${(d+\ell)^{n-1} \choose r}$ (resp. ${\ell \choose r}$) ways to form the pattern $\tau$. The rest is clear.
\end{proof}

The following example is a corollary to Theorem~\ref{multipattern}.

\begin{example}\label{ex4} If we consider the morphism $\phi_w$ defined in Section~\ref{introduction} and the pattern $\tau=1\mn 1\mn 1\mn 1$ then $d=\ell=1$, $r=4$, hence the number of occurrences of $\tau$ in $\phi^n(1)$ is 0, whenever $n=1$ or $n=2$, and $3\cdot{2^{n-1} \choose 4}$ otherwise.
\end{example}


\section{Patterns without internal dashes}

In what follows we need to extend the notion of an external occurrence of a pattern. Suppose $W=AXBYC$, where $A$, $X$, $B$, $Y$ and $C$ are some subwords. We say that an occurrence of $\tau$ in $W$ is external for a pair of words $(X,Y)$, if this occurrence starts in $X$, ends in $Y$ and is allowed to have some of its letters in $B$. For instance, if $W=12324245$, where $A=1$, $X=23$, $B=2$ and $Y=424$ then an occurrence of the generalized pattern $213$, namely the subword $324$ is an external occurrence for $(X,Y)$.

\begin{theorem}\label{generalTheorem_2}
Let $\mathcal{A}$ $=\{ 1 ,2,\ldots, k \}$ be an alphabet and a generalized pattern $\tau$ has no internal dashes. Let $X_1$ begins with the letter 1 and consists of $\ell$ copies of each letter $i\in \mathcal{A}$ {\rm(}$\ell \geq 1${\rm)}. Let a morphism $\phi$ be such that
$$1 \rightarrow X_1, \ 2 \rightarrow X_2, \ 3 \rightarrow X_3, \ldots , k \rightarrow X_k,$$
where we allow $X_i$ to be the empty word $\epsilon$ for $i=2,3,\ldots, k$ {\rm(}that is, $\phi$ may be an erasing morphism{\rm)}, $\displaystyle\sum_{i=2}^{k}|X_i|=k\cdot d$, and each letter from $\mathcal{A}$ appears in the word $X_2X_3\ldots X_k$ exactly $d$ times. Besides, we assume that there are no external occurrences of $\tau$ in $\phi^n(1)$ for the pair $(X_i,X_j)$ for each $i$ and $j$. Let $s_i$ be the number of internal occurrences of $\tau$ in $X_i$. In particular, $s_i=0$, whenever $X_i=\epsilon$. Then $N^{\tau}_{\phi}(1)=s_1$ and for $n\geq 2$, $N^{\tau}_{\phi}(n)=(d+\ell)^{n-2}\sum_{i=1}^{k}s_i$.
\end{theorem}

\begin{proof}
The theorem is straightforward to prove by observing that for $n\geq 2$, $\phi^n(1)$ has $(d+\ell)^{n-2}$ occurrences of each word $X_i$ (see the proof of Theorem~\ref{generalTheorem}).
\end{proof}

\begin{remark} In order to use Theorem~\ref{generalTheorem_2}, we need to control the absence of external occurrences of a pattern $\tau$ for given $\tau$ (without internal dashes) and a morphism~$\phi$. To do this, we need, for any pair $(X_i,X_j)$, to consider all the words $X_iWX_j$, where $|W|<|\tau|-1$, and $W$ is a permutation of a number of words from the set $\{ X_1, X_2, \ldots, X _k \}$. \end{remark}

The following examples are corollaries to Theorem~\ref{generalTheorem_2}.

\begin{example}\label{ex5} If we consider the morphism $\phi_w$ defined in Section~\ref{introduction} and the pattern $\tau=12$ then all the conditions of Theorems~\ref{generalTheorem_2} hold. In this case $d=\ell=1$ and $s=(2,1,0)$. Hence, the number of occurrences of the patterns $12$, that is the number of rises, is given by $N^{12}_{\phi_w}(1)=2$ and, for $n\geq 2$, $N^{12}_{\phi_w}(n)=3\cdot 2^{n-2}$. If $\tau=123$ then we can apply the theorem to get that for $n\geq 2$, $N^{123}_{\phi_w}(n)=2^{n-2}$.

If we want to count the number of occurrences of the pattern $\tau=21$, that is the number of descents, then we cannot apply Theorem~\ref{generalTheorem}, since for instance, the pair $(X_1, X_2) = (123, 13)$ has an external occurrence of $\tau$. However, it is obvious that the number of descents in $\phi^n(1)$ is equal to $|\phi^n(1)|-N^{12}_{\phi_w}(1)-1=3\cdot 2^{n-2}-1$.
\end{example}

\begin{example} If we consider the morphism $\phi$: $1\rightarrow 1243$, $2\rightarrow 3$, $3\rightarrow \epsilon$, and $4\rightarrow 124$ then for the pattern $\tau=123$, all the conditions of Theorems~\ref{generalTheorem_2} hold. In this case $d=\ell=1$, $s=(1,0,0,1)$. Hence, for $n\geq 1$, $N^{123}_{\phi}(n)=2^{n-1}$. For $\tau=321$ we cannot apply Theorem~\ref{generalTheorem_2}, since the pair $(X_4,X_1)$ has an external occurrence of $\tau$ (look at $X_4X_2X_1=12{\bf 431}243$). Consideration of the words $X_4X_2$ and $X_4X_1$ implies that the theorem cannot be apply for the patterns $132$ and $231$ respectively. However, we can apply the theorem to the pattern $213$ to prove that it does not occur in $\phi^n(1)$ for any~$n$.
\end{example}

{\bf Acknowledgement:} The final version of this paper was written
during the second author's (T.M.) stay at Haifa University, Haifa
31905, Israel. T.M. wants to express his gratitude to Haifa
University for the support.

\end{document}